\documentclass[12pt]{article}
\usepackage{amsfonts}
\usepackage{amsmath}
\usepackage{epsfig}

\title{On Nearly Optimal Paper Moebius Bands}
\author{Richard Evan Schwartz \thanks{Supported by N.S.F. Grant DMS-2102802 and also a Mercator Fellowship.}}

\newtheorem{theorem}{Theorem}[section]

\newtheorem{lemma}[theorem]{Lemma}

\newtheorem{corollary}[theorem]{Corollary}

\def\startproof{{\bf {\medskip}{\noindent}Proof: }}

\def\endproof{$\spadesuit$  \newline}

\def\R{\mbox{\boldmath{$R$}}}%
\def\Z{\mbox{\boldmath{$Z$}}}%

\begin{document}
\maketitle

\begin{abstract}
  Let $\epsilon<1/384$ and let
  $\Omega$ be a smooth embedded paper Moebius
  band of aspect
  ratio less than $\sqrt 3 + \epsilon$.  We prove that
  $\Omega$ is within Hausdorff distance
  $18 \sqrt \epsilon$ of an equilateral triangle of
  perimeter $2 \sqrt 3$.  This 
  is an effective and fairly sharp version of our recent theorems in [{\bf S0\/}] about
  the optimal paper Moebius band.
        \end{abstract}

\section{Introduction}

A {\it smooth embedded paper Moebius band\/} is an
infinitely differentiable isometric embedding
$I: M_{\lambda} \to \R^3$, where
$M_{\lambda}$ is the flat Mobius band
obtained by identifying the length-$1$ sides of a
$1 \times \lambda$ rectangle. 
We set $\Omega=I(M_{\lambda})$.  The number
$\lambda$ is the {\it aspect ratio\/} of $\Omega$.
This terminology is a mouthful, so I will use
the shorter term {\it paper Moebius band\/} to mean a smooth
embedded paper Moebius band.

In [{\bf S0\/}] I proved that
 a paper Moebius band has aspect ratio
 greater than $\sqrt 3$.  This bound was conjectured in 1977 by
 Halpern and Weaver in [{\bf HW\/}].  The article
 [{\bf FT\/}, \S 14] discusses the history of this conjecture and
 gives context for it.  [{\bf S0\/}] also has a discussion with
 many references.
 
 I also proved a limiting result in [{\bf S0\/}]:
   Let $\{I_n: M_{\lambda_n} \to \Omega_n\}$ be a sequence paper
   Moebius bands such that $\lambda_n \to \sqrt 3$.
   Up to isometry, $I_n$ converges
   uniformly to the triangular Moebius band map
   shown in Figure 1.
   \begin{center}
\resizebox{!}{1in}{\includegraphics{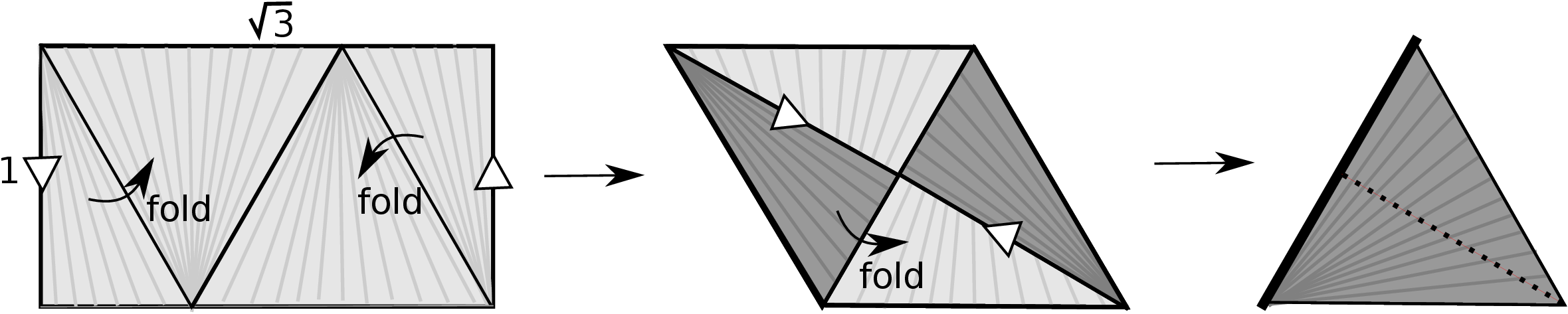}}
\newline
{\bf Figure 1:\/} The triangular Moebius band.
\end{center}

The purpose of this article is to make the results from [{\bf S0\/}]
more  effective.
Let $I_0: M_{\sqrt 3} \to \Omega_{\sqrt 3}$ be the
triangular Moebius band.
Given some other map $f: M_{\sqrt 3} \to \R^3$ we let
\begin{equation}
  \label{boundary}
  \|I_0-f\|_{\infty}=\sup_{p \in \partial M_{\sqrt 3}}
  \|I_0(p)-f(p)\|.
\end{equation}
Note that in Equation \ref{boundary} we are only
taking the sup over the boundary.

\begin{theorem}
  \label{eff}
  If $\epsilon<1/4$ and
  $I: M_{\lambda}\to \Omega$ is a
  paper Moebius band with
  $\lambda<\sqrt 3 + \epsilon$ then
  there is a homeomorphism
  $\phi: \partial M_{\sqrt 3} \to \partial M_{\lambda}$ and an
  isometry $\psi: \R^3 \to \R^3$ such that
  $\|I_0-\psi \circ I \circ \phi\|_{\infty}<6\sqrt \epsilon$.
\end{theorem}

Theorem \ref{eff} only deals with the boundary
of $\Omega$.  Here is another result, with somewhat
more restrictive conditions, that deals with the whole
space.

\begin{theorem}
  \label{eff2}
Let  $I: M_{\lambda}\to \Omega$ be a
  paper Moebius band with
  $\lambda<\sqrt 3 + \epsilon$.  If $\epsilon<1/4$ then
  there is an equilateral triangle
  $\bigtriangledown_0$ of perimeter $2 \sqrt 3$
  such that every point of $\Omega$ is
  within  $6 \sqrt \epsilon$ of $\bigtriangledown_0$.
  If $\epsilon<1/384$, then every point of
  $\bigtriangledown_0$ is within $18\sqrt \epsilon$
  of $\Omega$.
\end{theorem}

Let me phrase Theorem \ref{eff2} more gracefully.
The {\it Hausdorff distance\/} between two compact
subsets of $\R^n$ is defined to be the minimal
$\epsilon$ such that each set is contained in the
$\epsilon$-tubular neighborhood of the other.
(The minimum exists by compactness.)
The Hausdorff distance gives a well-known metric
on the set of compact subsets of $\R^n$.
We are working with $n=3$ in this paper.

\begin{corollary}[Main]
  Suppose $\epsilon<1/384$ and
  $I: M_{\lambda}\to \Omega$ is a
  paper Moebius band with
  $\lambda<\sqrt 3 + \epsilon$.  Then,
  in the Hausdorff metric,
  $\Omega$ is within $18 \sqrt \epsilon$ of an
  equilateral triangle of perimeter $2\sqrt 3$.
  \end{corollary}

\noindent
{\bf Remarks:\/}
(1) In Theorem \ref{eff}, I might have tried to get control
over the entire map and not just the boundary map.
I think it would be possible to do, 
with more effort and perhaps a worse estimate, but I wanted to keep
the proof light.
(2) The
cutoff of 1/4 is somewhat arbitrary, but it seems
like a reasonable notion
of ``fairly small''.  The significance of $\epsilon<1/384$ is
that then $6\sqrt{\epsilon}<1/3$, the in-radius of the
triangle $\bigtriangledown_0$.  The fact that $384=18^2$
is not very significant; it
is an artifact of the proofs.
(3) I don't know how sharp the constants in Theorems \ref{eff} and
\ref{eff2} are but the $O(\sqrt \epsilon)$ term is sharp.
This derives from the fact that a right triangle with
base $1$ and hypotenuse $1+\epsilon$ has height
$O(\sqrt \epsilon)$.  See \S \ref{ex} for more about this.
\newline

In \S 2 I will prove a few easy technical lemmas that will
help with the overall proof, and then in \S 3 I will prove
Theorems \ref{eff} and \ref{eff2}.
To make this paper more self-contained, I recall the main
arguments in [{\bf S0\/}] in some detail.
In \S \ref{ex} I give an example illustrating
the sharpness of our results.

In [{\bf S0\/}] I thanked many people for
their help and insights into paper Moebius bands.
I thank all these people again.
This research is supported by a Simons
Sabbatical Fellowship, a grant from
the National Science Foundation, and
a Mercator Fellowship. I'd like to
thank all these organizations.

      \section{Some Perturbation Results}

We let $\ell(\cdot)$ denote
arc length, and we assume
that $0<\epsilon<1/4$.

\subsection{Two Trivialities}

We will need two very easy estimates about the square root function.
We single these out in advance to make the rest of our exposition
go more smoothly. 

\begin{lemma}
  \label{sq0}
  If $L<3/2$ then
  $\sqrt{L^2+(13/4)\epsilon}>L+\epsilon$.
\end{lemma}

\startproof
$(L+\epsilon)^2=L^2+(2L+\epsilon)\epsilon<L^2+(13/4)\epsilon$.
\endproof

\begin{lemma}
  \label{sq1}
  If $L<3/\sqrt 2$ then
  $\sqrt{L^2+(9/2)\epsilon}>L+\epsilon/2$.
\end{lemma}

\startproof
$(L+(\epsilon/2))^2=L^2+(L+(\epsilon/4))\epsilon<L^2+(9/2)\epsilon.$
\endproof

\subsection{Perturbing an Isosceles Triangle}

Let $\bigtriangledown$ be a triangle with a horizontal base.
We say that the {\it bottom vertex\/} of $\bigtriangledown$ is
the vertex not on the base.
Let $\bigtriangledown_*$ denote the isosceles triangle having
the same base and height as $\bigtriangledown$.  We get
$\bigtriangledown_*$ from $\bigtriangledown$ by moving
the bottom vertex horizontally by a distance we call
$d(\bigtriangledown,\bigtriangledown_*)$.
Let $\vee$ denote the union of the non-horizontal edges of
$\bigtriangledown$.  Likewise define $\vee_*$.

\begin{lemma}
  \label{offset1}
  Suppose $\ell(\vee_*)<3$ and the slopes of the sides of
  $\vee_*$ exceed $1$ in absolute value.
  If $d(\bigtriangledown,\bigtriangledown_*) \geq
  \sqrt{13\epsilon/2}$ then
 $\ell(\vee)>\ell(\vee_*)+2\epsilon.$
\end{lemma}

\startproof
Let $p_1,p_2,q$ be the vertices of $\bigtriangledown$, with
$q$ being the bottom vertex.  Likewise let $q^*$ be the bottom
vertex of $\bigtriangledown_*$.  Let $r_2$ denote the reflection of
$p_2$ in the horizontal line through $q$ and $q^*$. 
  By symmetry, $\ell(\vee_*)=\ell(\overline{p_1r_2})$ and
  $\ell(\vee)=\ell(\overline{p_2q})+\ell(\overline{qr_2})$.
  By assumption, $\ell(\overline{p_1r_2})<3$ and
  the horizontal distance from $q$ to $\overline{p_1r_2}$ is at least $\sqrt{13\epsilon/2}$.
  Since $\overline{p_1r_2}$ has slope at least $1$ in absolute value, the distance
  from $q$ to $\overline{p_1r_2}$ is at least $d=\sqrt{13\epsilon/4}$.
  
\begin{center}
\resizebox{!}{1.75in}{\includegraphics{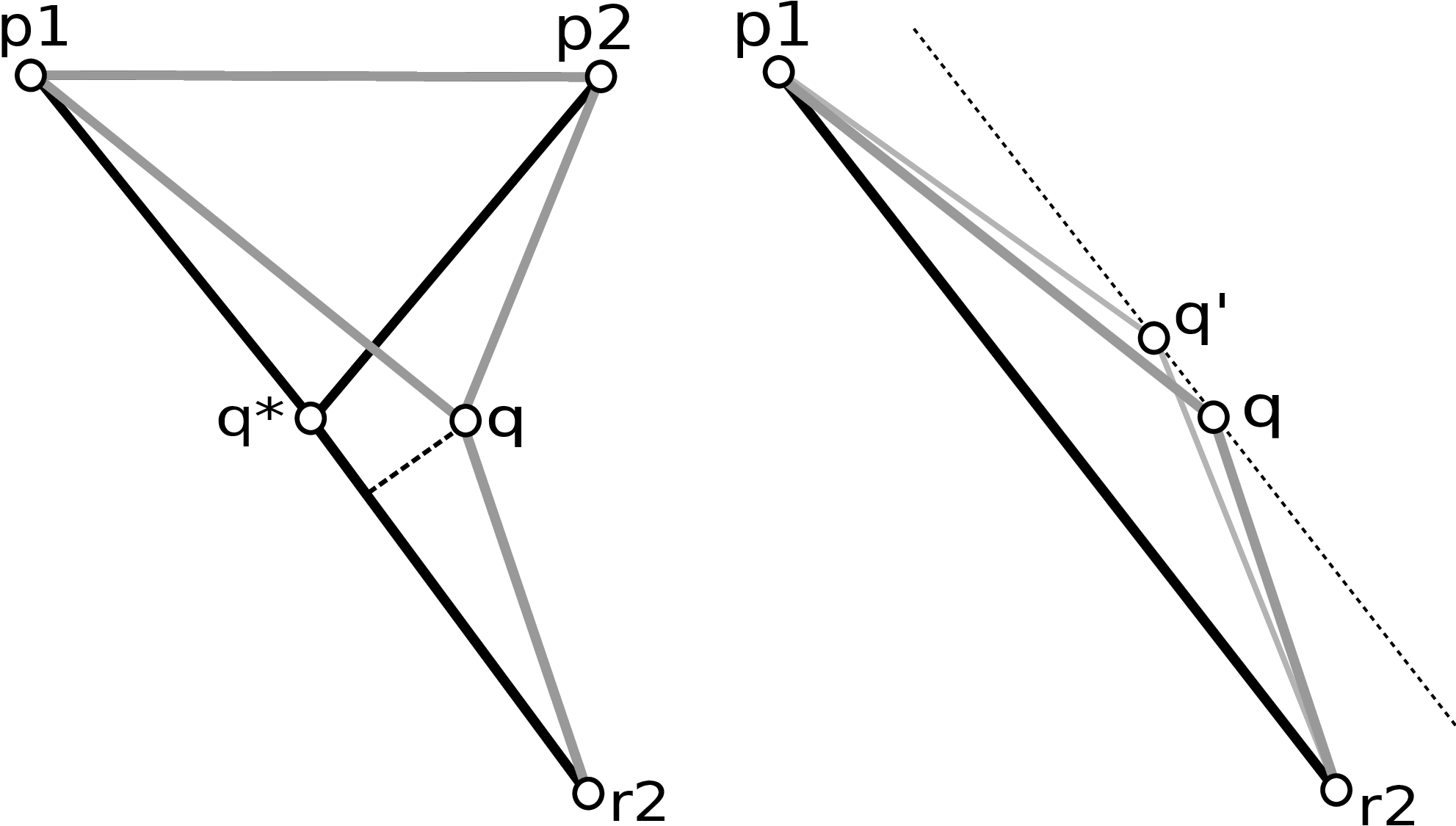}}
\newline
{\bf Figure 2:\/}  The quantities used in the proof.
  \end{center}

  Let $L=(1/2)\ell(\vee_*)<3/2$.
 We form a new triangle $T$ by sliding $q$ parallel to $\overline{p_1r_2}$,
  to a new point $q'$ which minimizes
  $\ell(\overline{p_1q'})+\ell(\overline{q'r_2})$.
  This happens when the triangle $(p_1,r_2,q')$ is isosceles.
  We have
  $$\ell(\overline{p_1q'})=\ell(\overline{q'r_2}) \geq
  \sqrt{L^2+d^2}=\sqrt{L^2+(13/4)\epsilon}>L+\epsilon.$$
The last inequality is Lemma \ref{sq0}.
Hence
  $$\ell(\vee)=\ell(\overline{p_1q})+\ell(\overline{qr_2}) \geq 
  \ell(\overline{p_1q'})+\ell(\overline{q'r_2})>2L+2\epsilon=\ell(\vee_*)+2\epsilon.$$
This completes the proof.
  \endproof

  \noindent
  {\bf Remark:\/}
  We call the trick of replacing $q$ by $q'$ the {\it sliding
    trick\/}.
  We will use this trick in the next section as well.

\subsection{Perturbing an Affine Map}

Let $S$ be a line segment with $\ell(S)<3$.

   Let $I: S \to \R^3$ be a
   unit speed map.  Let $I_*: S \to \R^3$ be the
   affine map that agrees with $I$ on the endpoints.
   Let     $S'=I(S)$ and
   $ S'_*=I_*(S)$.
      Note that $\ell(S'_*) \leq \ell(S')=\ell(S)<3$.

  Consider the graph $\Gamma$ of $I$ defined as
   $$\Gamma=\bigcup_{x \in S} (x,I(x)) \subset \R^4= \R \times \R^3.$$
   Likewise define $\Gamma_*$.  The graph $\Gamma_*$ is
   a straight line segment.

   \begin{lemma}
     \label{graph}
     $\ell(\Gamma_*) \leq \ell(\Gamma) \leq 3 \sqrt 2$ and 
        $\ell(\Gamma)-\ell(\Gamma_*) \leq   \ell(S')-\ell(S'_*)$.
   \end{lemma}

   \startproof
   The maps $I$ and $I_*$ respectively have speeds $1$ and
   $s=\ell(S'_*)/\ell(S')$.
   As curves,  $\Gamma_*$ and $\Gamma$ respectively have speeds
   $\sqrt{1+s^2}$ and $\sqrt 2$.
   Integrating these functions over the domain, which has length $\ell(S')$, we
   see that
   \begin{equation}
     \ell(\Gamma_*)=\ell(S') \sqrt{1+s^2}=\sqrt{\ell(S')^2+\ell(S'_*)^2}<3 \sqrt 2.
   \end{equation}
   \begin{equation}
     \ell(\Gamma)=\ell(S') \sqrt 2=\sqrt{\ell(S')^2+\ell(S')^2}<3
     \sqrt 2.
   \end{equation}
   That takes care of the first statement.

   For the second statement,
   let $a=\ell(S')^2$ and let
$f(t)=\sqrt{a+t^2}.$
We have $|f'(t)| \leq 1$ for all $t$, regardless of the value of $a
\geq 0$.
This means that for $0 \leq t_*<t$ we have $f(t)-f(t_*) \leq t-t_*$.
Applying this to $t=\ell(S')$ and $t_*=\ell(S'_*)$, we get
$\ell(\Gamma)-\ell(\Gamma_*)\leq \ell(S')-\ell(S'_*)$.
\endproof

Let $\|\cdot \|_{\infty}$ denote the sup-norm on $S$.
      
   \begin{lemma}
     \label{wiggle}
     If $\|I-I_*\|_{\infty} \geq 3 \sqrt{\epsilon}$
      then $\ell(S')>\ell(S_*')+\epsilon$.
   \end{lemma}

   \startproof
   By hypothesis, there is some
   $p \in S$ such that $\|I_*(p)-I(p)\| \geq 3\sqrt{\epsilon}$. Hence,
   there is a
   point $q \in \Gamma$ which is at least $3\sqrt{\epsilon}$
   away from the point in $\Gamma_*$ having the same
   $\R$ coordinate.    Since $I_*$ has speed $s<1$, the graph
   $\Gamma_*$ has ``slope'' at most $1$.
   This implies that
   $q$ is at least $\sqrt{9 \epsilon/2}$ from
   $\Gamma_*$ in $\R^4$.

   Let $L=\ell(\Gamma_*)/2<3/\sqrt 2$.
   By the same kind of point-sliding trick we used in the proof
   of Lemma \ref{offset1}, we
  have
  \begin{equation}
    \ell(S')-\ell(S'_*) \geq
  \ell(\Gamma) -\ell(\Gamma_*) \geq
  2\sqrt{L^2+(9/2)\epsilon} -2L \geq \epsilon.
\end{equation}
The first inequality if Lemma \ref{graph}.
The last inequality is Lemma \ref{sq1}.
\endproof

    \section{Proofs of the Results}

\subsection{The Optimality Proof Revisited}

I first will recall some of the material in
[{\bf S0\/}] in order to give more self-contained
proof of Theorems \ref{eff} and \ref{eff2}.
Let $I: M_{\lambda} \to \Omega$
be a (smooth embedded) paper Moebius band.
A {\it bend\/} on $\Omega$ is a straight line segment in $\Omega$
which has its endpoints (and only its endpoints)
in $\partial \Omega$.

\begin{lemma}
  $\Omega$ has a foliation by bends.
\end{lemma}

\startproof
(Sketch)
This seems to be a folklore result.
See  [{\bf CL\/}],  [{\bf HM\/}], [{\bf Mas\/}] for
arguments which immediately work in case
$\Omega$ has an open dense set of points
with nonzero mean curvature.
See [{\bf S0\/}, Prop. 2.1] for the general case
and more precise references.
The basic idea is that each point on $\Omega$
of nonzero curvature has a unique tangent direction
where the differential of the Gauss map is trivial.
These directions integrate up to give a foliation
of the nonzero-mean-curvature subset of
$\Omega$ by line segments.  The
complementary pieces are planar trapezoids
and one can
fill in this foliation in an obvious way.
\endproof

A (embedded) $T$-{\it pattern\/} is a
pair of disjoint bends on the paper Moebius
band which lie in perpendicular intersecting lines.

\begin{lemma}
  $\Omega$ has a $T$-pattern.
  \end{lemma}

  \startproof
  (Sketch)
  This is [{\bf S0\/}, Lemma T].
  The space of bends in the bend
  foliation is homeomorphic to a circle.
  The space of pairs of unequal bends
  has a $2$-point compactification which
  makes it homeomorphic to the sphere $S^2$.
  Using the dot product and the cross product
  we define $2$ odd functions on $S^2$ which,
  when they simultaneosly vanish, detect a
  $T$-pattern.   By the Borsuk-Ulam Theorem,
  both functions do simultaneously vanish,
  and this gives the $T$-pattern.
  \endproof
  
[{\bf S0\/}, Lemma G] uses the $T$-pattern to show that
$\lambda>\sqrt 3$.   We essentially reproduce this argument here.

\begin{center}
\resizebox{!}{1.7in}{\includegraphics{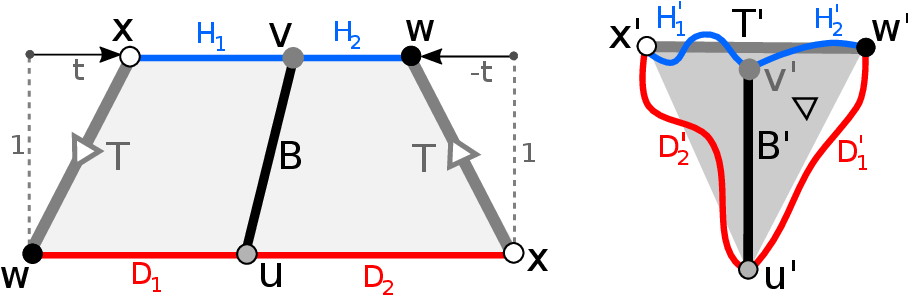}}
\newline
    {\bf Figure 3:\/} The trapezoid $\tau$ (left) and the T-pattern (right).
  \end{center}

Let $S'=I(S)$ for any $S \subset M_{\lambda}$.
We rotate $\Omega$ so that one of the bends of the
$T$-pattern, $T'$, lies in
$X$-axis and the other bend, $B'$, lies in the negative ray of the
$Y$-axis.  Let
$B=I^{-1}(B')$ and $T=I^{-1}(T')$.
These preimages are also line segments.
We cut $M_{\lambda}$ open along the line segment
$T$ to get a bilaterally
symmetric trapezoid $\tau$.   We normalize $\tau$
so that the parallel
sides  are horizontal and so that
$u,v,w,x$ are
mapped to $\Omega$ as in Figure 3.  
The quantity $t$ is the horizontal
displacement of $T$.  Figure 3 shows the
case when $t>0$.

Note that the shaded triangle $\bigtriangledown$ on
the right side of Figure 3 has
base $\ell(T')=\sqrt{1+t^2}$ and height greater $\ell(B) \geq 1$.
Let $\vee$ denote the union of the two non-horizontal
sides of $\bigtriangledown$.  The easy   [{\bf S0\/}, Lemma 2.2]
says that $\ell(\vee)>\sqrt{5+t^2}$.

  We have the following equations
$$
    \lambda=\ell(H)+t=\ell(D)-t,
    $$
    $$
    \ell(H)=\ell(H')>\ell(T')=\ell(T)=\sqrt{1+t^2}
$$
\begin{equation}
  \label{key}
    \ell(D)=\ell(D') \geq \ell(\vee)>\sqrt{5+t^2}.
\end{equation}

 By Equation \ref{key},
 \begin{equation}
   \label{bigmax}
   \lambda> \max(h(t),d(t)), \hskip 20 pt
   h(t)=\sqrt{1+t^2}+t, \hskip 20 pt
   d(t)=\sqrt{5+t^2}-t.
 \end{equation}
 Taking derivatives, we see easily that $h$ is an increasing
 function and $d$ is a decreasing function.  Also
 $h(1/\sqrt 3)=d(1/\sqrt 3)=\sqrt 3$.  All this implies that
 $h(t) \geq \sqrt 3$ if $t \geq 1/\sqrt 3$ and
 $d(t) \geq \sqrt 3$ if $t \leq 1/\sqrt 3$.
 Hence $\lambda>\sqrt 3$ regardless of the value of $t$.
 That proves Lemma G.

 \subsection{Geometric Bounds}
 
 Now we go beyond what we did in [{\bf S0\/}].
   If $S$ is any object associated to our paper Moebius
  band $\Omega$, let $S_0$ be the corresponding object associated to
  the triangular Moebius band $\Omega_0$.  In particular,
  let $(B'_0,T'_0)$ be the $T$-pattern associated to $\Omega_0$.
  These are normalized as above.  We take $\psi$ to be the
  translation so that $\psi(T')$ and $T'_0$ have the same midpoint.
  (Both are segments in the $X$-axis.)
  Let $t_0=1/\sqrt 3$.

 \begin{lemma}
   \label{lip}
   $|t-t_0|<4 \epsilon/3$.  In particular $t \in (0,1)$.
 \end{lemma}

 \startproof
The functions $h$ and $d$ are
both convex, and $h'(t_0)=3/2$ and $d'(t_0)=-3/4$.
By convexity,
 $|h'(t)| \geq 3/2$ when $t \geq t_0$ and
 $|d'(t)| \geq 3/4$ when $t \leq t_0$.
 Under our assumption that $\lambda<\sqrt 3+\epsilon$ we
 get $|t-t_0|<4 \epsilon/3 \leq 1/3$.  The second bound now
 follows from the fact that $1/3<t_0<2/3$.
 \endproof

 \begin{lemma}[Length Bound]
     $\ell(H)<\ell(D)<3$.
   \end{lemma}
 
\startproof
  Since $t>0$ we have
   $\ell(H)<\ell(D)$.   
   Since $\epsilon \leq 1/4$ we have $\lambda<2$.
   Since $t<1$ have
   $\ell(D)=\lambda+t<2+1=3$.
      \endproof

    \begin{lemma}
      \label{base}
      $|\ell(T')-2t_0|<\epsilon$.
    \end{lemma}

    \startproof
    We have $\ell(T')=f(t):=\sqrt{1+t^2}$ and
    $|t-t_0|<4 \epsilon/3$.  Also, we have $f(t_0)=2t_0$.
    At the same time, $|f'(s)|<3/4$ on $(0,1)$, an interval that
    contains $t$ and $t_0$.  Hence,  Lemma \ref{lip} gives
    $|\ell(T')-2t_0|<\epsilon$, as claimed.
    \endproof

    \begin{lemma}
      \label{height}
      The height $y$ of $\bigtriangledown$ is less than
      $1+\epsilon$.
    \end{lemma}

    \startproof
    Let  $d(y,t)=\sqrt{1+t^2+4y^2}-t$.
    For each choice of $y$ the function $d(*,t)$ is  decreasing.
    The obvious generalization of our proof of Lemma G above gives
    $\lambda>h(t_y)$ where $t_y$ is such that
    $h(t_y)=d(y,t_y)$.  We have $t_y=\pm y^2/\sqrt{1+2y^2}$.
      The negative choice leads to $\lambda>\sqrt 5>\sqrt 3 +
      \epsilon$.
     The positive choice leads to
      $\lambda>g(y):=h(t_y)=\sqrt{1+2y^2}$.
      Now $g(1)=\sqrt 3$ and $g'(y)>1$ for $y \geq 1$.
      Since $\lambda<\sqrt 3 + \epsilon$ we have
      $y<1+\epsilon$.
\endproof

Let $\bigtriangledown_*$be as in
Lemma \ref{offset1}.  Let
$\delta=d(\bigtriangledown,\bigtriangledown_*)$.

\begin{lemma}
  \label{offset}
  $\delta<\sqrt{13\epsilon/2}$.
\end{lemma}

\startproof
We use the notation from Lemma \ref{offset1}.
Assume that $\delta \geq  \sqrt{13\epsilon/2}$.
Note that $$\ell(\vee_*) \leq \ell(\vee) \leq \ell(D') =\ell(D)<3.$$
The base of $\bigtriangledown$ is $\sqrt{1+t^2}<2$, and
the height of $\bigtriangledown$
is at least $1$. Hence the absolute values of the
slopes of the sides of $\vee_*$ are at least $1$.
Since $\delta \geq \sqrt{13 \epsilon/2}$ we now have
all the hypotheses of Lemma \ref{offset1}, which tells us that
$\ell(\vee)-\ell(\vee_*)>2\epsilon$.  This also tells us that
$\ell(\bigtriangledown)-\ell(\bigtriangledown_*)>2\epsilon$.

If $t \geq t_0$ then $\ell(\bigtriangledown_*) \geq 2 \sqrt 3$, and so
$$2\lambda=\ell(H)+\ell(D)=\ell(H')+\ell(D')
  \geq
  \ell(\bigtriangledown)>\ell(\bigtriangledown_*)+2\epsilon
  \geq 2 \sqrt 3 + 2\epsilon.$$
If $t<t_0$ then we consider the function $d(t)$ from Equation
\ref{bigmax} and observe that
  $$\lambda \geq \ell(\vee)-t >
  \ell(\vee_*)-t+2\epsilon>
  d(t)+2\epsilon  \geq \sqrt 3+  2 \epsilon.$$
  Either case gives $\lambda>\sqrt 3+\epsilon$, a contradiction.
Hence    $\delta<\sqrt{13\epsilon/2}$.
  \endproof

  \begin{lemma}[T Pattern Bound]
    Each endpoint of the $T$-pattern $(T,B)$ is within
    $3 \sqrt \epsilon$ of the corresponding endpoint of
    $(T_0,B_0)$. \end{lemma}

 \startproof
 By Lemma \ref{base}, we have
 $\|x'_0-x'\|<\epsilon/2$ and
 $\|w'_0-w'\|<\epsilon/2$.
 This takes care of the endpints of
 $T'$.

 Since $\bigtriangledown$ has height less than
 $1+\epsilon$, the $y$-coordinates of the
 endpoints of $B'$ are
 within $\epsilon$ of the
 $y$-coordinates of the corresponding endpoints of
 $B'_0$.   From Lemma \ref{offset}, the
 $x$-coordinates of the
 endpoints of $B'$ are within
 $\sqrt{13 \epsilon/2}$ of the 
 $x$-coordinates of $B_0'$.
 By the Pythagorean Theorem,
   \begin{equation}
   \|u'-u'_0\|,  \|v'-v_0'\|<\sqrt{(13/2)\epsilon+\epsilon^2}<3 \sqrt \epsilon.
  \end{equation}
  The final inequality is equivalent to the statement that
  $p(\epsilon)=(5/2)\epsilon - \epsilon^2$ is positive
  for $\epsilon \in (0,1/4)$, and this is certainly true.
 \endproof

 \subsection{Proof of Theorem \ref{eff}}

 We have our paper Moebius band
 $$I: M_{\lambda} \to \Omega.$$
   We have quadrilaterals $\tau$ and $\tau_0$ as
   in the left side of Figure 3.  We think of $\partial \tau$ as
   having $6$ edges:  In addition to
   $H_1,H_2,D_1,D_2$ we have the non-horizontal
   edges which we call $T_1$ and $T_2$.
   By construction $T_1'=T_2'=T'$.
   We make the same labelings for $\tau_0$.
We define   $I_*: \partial \tau \to \R^3$ to be the unique
map that is affine on each edge of $\tau$.
With this definition,  $I_*$ induces a piecewise linear
map on the union $\partial M_{\lambda}$.

   \begin{lemma}
     \label{linebound}
     $\|I_*-I\|_{\infty}<3\sqrt \epsilon$.
  \end{lemma}

  \startproof
  Here, as in Theorem \ref{eff}, we are taking the sup-norm over points
  on $\partial M_{\lambda}$.
     For each
   $S \in \{H_1,H_2,D_1,D_2\}$ we have $\ell(S)<3$.
   Since $\lambda<\sqrt 3+\epsilon$ we must have
   $\ell(S')<\ell(S'_*)+\epsilon$.  This is to say that
   we cannot have more than $\epsilon$ of slack in
   any of these curves.
   Now we conclude, by
   Lemma \ref{wiggle}, that
   $$\sup_{p \in S} \|I_*(p)\|<3 \sqrt \epsilon.$$
   Since this works for all choices of $S$, we get the bound
   of this lemma.
   \endproof

  We define $\phi: \partial M_{\sqrt 3} \to \partial M_{\lambda}$ to be the
  piecewise linear map which maps each edge of
  $\tau_0$ to the corresponding edge of $\tau$.
  Lemma \ref{linebound} gives us
    
  \begin{equation}
    \label{est1}
  \|I \circ \phi - I_* \circ \phi\|_{\infty}=
  \|I-I_*\|_{\infty}<         3 \sqrt{\epsilon}.
   \end{equation}

   Finally, let us compare $I^* \circ \phi$ to $I_0$.  Both these
   maps are affine on each of the edges of $\tau_0$.
   Since these maps differ by at most $3\sqrt{\epsilon}$ on the
   endpoints of each edge of $\tau$, we have
   $\|I_0- I_* \circ \phi\|_{\infty} < 3 \sqrt \epsilon$.
   Combining this with Equation \ref{est1} we have
    \begin{equation}
   \|I_0- I \circ \phi\|_{\infty}  \leq 
   \|I_0- I_*\circ \phi\|_{\infty}  +
   \|I_* \circ \phi- I \circ \phi\|_{\infty}  <
   (3 + 3) \sqrt \epsilon<
   6 \sqrt \epsilon.
 \end{equation}
 This completes the proof of Theorem \ref{eff},

 \subsection{Proof of Theorem \ref{eff2}}

By Theorem \ref{eff}, each of these
bends foliating $\Omega$ has its endpoints
inside the $6 \sqrt{\epsilon}$-tubular neighborhood $N$ of 
the  triangle $\bigtriangledown_0$.
But $N$ is convex, and hence all the
bends foliating $\Omega$ lie within $N$.   Hence
$\Omega \subset N$.
This proves that every point of $\Omega$ is within
$6 \sqrt \epsilon$ of $\bigtriangledown_0$.
This is the first statement of Theorem \ref{eff2}.

For the second statement, we take $\epsilon<1/384$.  Note that
$d:= 6 \sqrt \epsilon < 1/3$,
and $1/3$ is the in-radius
of $\bigtriangledown_0$.
Figure 4 shows a picture of an annular neighborhood $A$ of
$\bigtriangledown_0$.   The core curve of
$A$ is $\bigtriangledown_0$.   
The disks, centered at the vertices of $\bigtriangledown_0$, each
have radius $d$. The fact that
$d<1/3$ guarantees that Figure 4 is a topologically
accurate picture.  Let $Dx'$ be the disk which is
centered at $x'$, etc.
Let $C$ be the solid triangle bounded
by the inner boundary of $A$.

\begin{center}
\resizebox{!}{1.9in}{\includegraphics{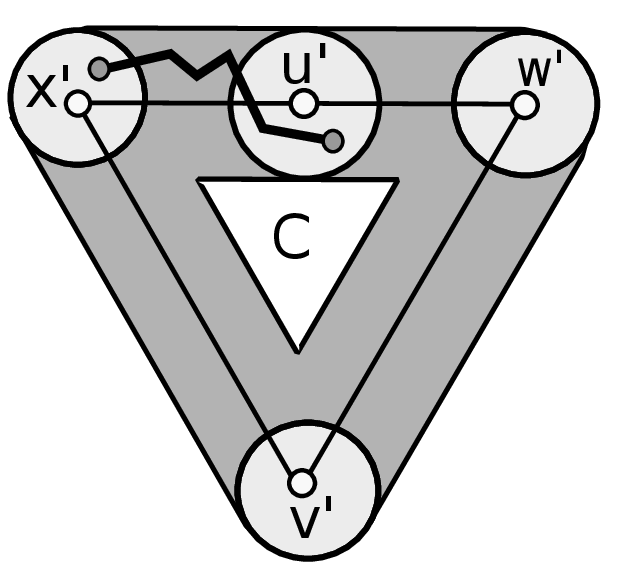}}
\newline
{\bf Figure 4:\/} The neighborhood $A$ of $\bigtriangledown_0$.
\end{center}

Let $f: \R^3 \to \R^2$ denote projection.
By Theorem \ref{eff}, we see that $f(\partial \Omega) \subset A$.
In particular, it makes to ask about which element in the
fundamental group $\pi_1(\R^2-C)$ this loop represents.

\begin{lemma}
  $f(\partial \Omega)$ generates $\pi_1(\R^2-C)$.
\end{lemma}

\startproof
  Consider the segment $H_1$ of $\partial \tau$.  The image $H_1'$
  is a curve which connects a point in $Dx'$ to a point in $Du'$ and
  remains in the portion of $A$ connecting these disks in the shortest
  way.   See Figure 4.
  In other words, $H'_1$ cannot ``wind the other way'' around $C$.
  The curves $H_2',D_1',D_2'$ have similar properties. From this we see
  that $f(\partial \Omega)$ and $f(\partial \Omega_0)$ represent the
  same
  element in  the fundamental group $\pi_1(\R^2-C)$, namely a
  generator.
    \endproof
  
  \begin{lemma}
    $C \subset f(\Omega)$.
  \end{lemma}

  \startproof
  Suppose not. Then there is
  $p \in C$ such that
  $f(\Omega) \subset \R^2-p$. 
  Note that $f(\partial \Omega)$ also generates
  $\pi_1(\R^2-p)$. 
  Let $\Lambda$ be the core curve
  of $\Omega$, the curve connecting
  all the midpoints of the bends. Note also that
  there is a homotopy from $f(\partial \Omega)$
  to a curve which winds twice around $f(\Lambda)$.
  We just push in along the bends.  Call
  this double-wrapped curve $f(2\Lambda)$.
  Since our homotopy avoids $p$, we see that
  $f(2\Lambda)$ generates $\pi_1(\R^2-p)=\Z$.
  This is impossible because
  $f(2\Lambda)$ is an even number in $\pi_1(\R^2-p)$.
  \endproof
  
  By the previous result, every vertical
    line through a point of $C$ intersects $\Omega$. But all
  of $\Omega$ lies in the $6 \sqrt \epsilon$ neighborhood of
  $\R^2$.  Hence every point of $C$ is within $6 \sqrt \epsilon$
  of a point of $\Omega$.  But every point of $\bigtriangledown_0$
  is within $12 \sqrt \epsilon$ of a point of $C$.  Hence
  every point of $\bigtriangledown_0$ is within $18 \sqrt \epsilon$
  of a point of $\Omega$.

  This completes the proof of
  Theorem \ref{eff2}.

  \subsection{Sharpness of the Results}
  \label{ex}

  Now we give an example illustrating the sharpness of our results.
  We will describe a polygonal paper Moebius band.  Using
    the {\it pseudo-fold method\/} in [{\bf HW\/}] (which just
    amounts to smoothing out the corners) we can then approximate
this object  as close as we like by smooth embedded paper Moebius
bands.
We omit the details of the smoothing.

  We start with $M_{\sqrt 3}$ and we insert a vertical strip $B$ of
  width $O(\epsilon)$ as shown on the left in Figure 5. 
   \begin{center}
\resizebox{!}{1.21in}{\includegraphics{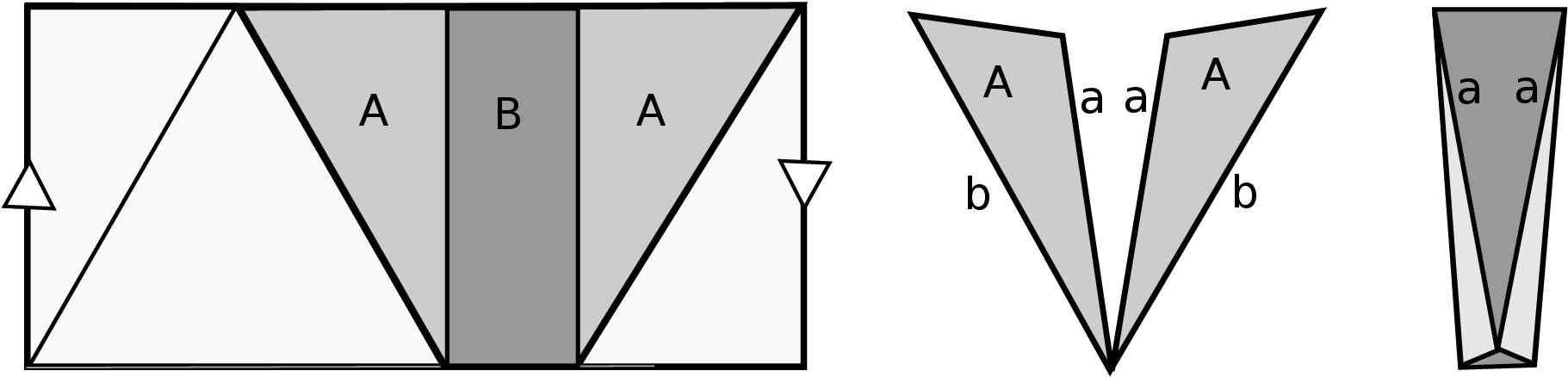}}
\newline
{\bf Figure 5:\/} Modified triangular Moebius band.
\end{center}
We start with the triangle $\bigtriangledown_0$ and
slit it open along its vertical midline.  This leaves us
with two triangular ``doors''.   Each door has a side
which coincides with a side of $\vee_0$.  Call these
the $b$-sides.  We rotate these
doors about the $b$-sides by an angle $O(\sqrt \epsilon)$.
The two vertical
sides come out of the plane and make a very sharp $V$.
We call this $V$ the {\it the crack\/}.
The width of the crack is $O(\epsilon)$ and
it rises $O(\sqrt \epsilon)$ out of the plane.
This derives from the fact that $1-\cos(\epsilon)=O(\epsilon^2)$
and $\sin(\epsilon)=O(\epsilon)$.
We have labeled the sides of the crack by the letter $a$.
We call the union of rotated triangles
{\it the saloon doors\/}.

At the same time, we fold up $B$ so that its two vertical edges
become the sides of the same kind of $V$.  We adjust the width,
keeping it $O(\epsilon)$, so that the new $V$ is isometric to the
original one, the one we call the crack.
 We call this folded image of $B$ (which is still planar)
{\it the plug\/}.
We now revolve the plug $\pi$ radians about the $Z$-axis, tilt
slightly, and glue it to the crack along the $a$-sides.
This creates a piecewise isometric map from
the $ABA$ quadrilateral on the left of Figure 5
to a region whose two outer bends comprise
$\vee_0$.   Our image, which we call
{\it the wrinkle\/}, rises $O(\sqrt \epsilon)$ out of the plane.

Finally, we
build the triangular Moebius band and replace the relevant
copy of $\bigtriangledown_0$ by the wrinkle.  For
the purposes of smooth approximation, we take care that
(in an infinitesimal sense) the wrinkle is ``on the outside'' so that
nearby approximations will be embedded.
By construction,
our new $\Omega$ has aspect ratio $\lambda+O(\epsilon)$,
contains $\bigtriangledown_0$, and also contains
points $O(\sqrt \epsilon)$ from $\bigtriangledown_0$.
The smoothings of $\Omega$ are the examples showing
the sharpness of $O(\sqrt \epsilon)$ in our results.

\section{References}

\noindent
  [{\bf CL\/}], S.-S. Chern and R. K. Lashof,
    {\it On the total curvature of immersed manifolds\/},
    Amer. J. Math. {\bf 79\/} (1957) pp 306--318
    \vskip 8 pt
        \noindent
[{\bf FT\/}], D. Fuchs, S. Tabachnikov, {\it Mathematical Omnibus: Thirty Lectures on Classic Mathematics\/}, AMS 2007
\vskip 8 pt
\noindent
 [{\bf HN\/}], P. Hartman and L. Nirenberg,
  {\it On spherical maps whose Jacobians do not change sign\/},
    Amer. J. Math. {\bf 81\/} (1959) pp 901--920
    \vskip 8 pt
    \noindent
[{\bf HW\/}], B. Halpern and C. Weaver,
{\it Inverting a cylinder through isometric immersions and embeddings\/},
Trans. Am. Math. Soc {\bf 230\/}, pp 41--70 (1977)
\vskip 8 pt
\noindent
[{\bf Mas\/}] W. S. Massey, {\it Surfaces of Gaussian Curvature Zero
  in Euclidean $3$-Space\/},  Tohoku Math J. (2) 14 (1), pp 73-79
(1962)
\vskip 8 pt
\noindent
[{\bf S0\/}] R. E. Schwartz, {\it The optimal paper Moebius band\/},
Annals of Math. (2024) to appear (see also arXiv 2308:12641)

\end{document}